\documentclass{amsart}
\usepackage{amscd, amssymb, amsfonts}
\usepackage{amsmath,amsthm}
\usepackage{amssymb}

  [2002/01/19 v2.2g %
    AMS font definitions%
  ]
\DeclareFontFamily{U}{euf}{}
\DeclareFontShape{U}{euf}{m}{n}{%
  <5><6><7><8><9>gen*eufm%
  <10><10.95><12><14.4><17.28><20.74><24.88>eufm10%
  }{}
\DeclareFontShape{U}{euf}{b}{n}{%
  <5><6><7><8><9>gen*eufb%
  <10><10.95><12><14.4><17.28><20.74><24.88>eufb10%
  }{}

  [2002/01/19 v2.2g %
    AMS font definitions%
  ]
\DeclareFontFamily{U}{msb}{}
\DeclareFontShape{U}{msb}{m}{n}{%
  <5><6><7><8><9>gen*msbm%
  <10><10.95><12><14.4><17.28><20.74><24.88>msbm10%
  }{}

  [2002/01/19 v2.2g %
    AMS font definitions%
  ]
\DeclareFontFamily{U}{msa}{}
\DeclareFontShape{U}{msa}{m}{n}{%
  <5><6><7><8><9>gen*msam%
  <10><10.95><12><14.4><17.28><20.74><24.88>msam10%
  }{}

\newtheorem{theorem}{Theorem}[section]

\newtheorem{proposition}[theorem]{Proposition}
\newtheorem{corollary}[theorem]{Corollary}

\theoremstyle{definition}

\theoremstyle{remark}

\numberwithin{equation}{section}

\begin{document}

\title[]
{Some identities of Bernoulli numbers and polynomials associated with Bernstein polynomials}

\author{Min-Soo Kim, Taekyun Kim, Byungje Lee and Cheon-Seoung Ryoo}

\begin{abstract}
We investigate some interesting properties of the Bernstein polynomials related to the
bosonic $p$-adic integrals on $\mathbb Z_p.$
\end{abstract}

\address{Department of Mathematics, KAIST, 373-1 Guseong-dong, Yuseong-gu, Daejeon 305-701, South Korea}
\email{minsookim@kaist.ac.kr}

\address{Division of General Education-Mathematics, Kwangwoon University, Seoul, 139-701, South Korea}
\email{tkkim@kw.ac.kr}

\address{Department of Wireless Communications Engineering, Kwangwoon University, Seoul 139-701, South Korea}
\email{bj\_lee@kw.ac.kr}

\address{Department of Mathematics, Hannam University, Daejeon 306-791, South Korea}
\email{ryoocs@hnu.kr}



\maketitle



\def\ord{\text{ord}_p}
\def\C{\mathbb C_p}
\def\BZ{\mathbb Z}
\def\Z{\mathbb Z_p}
\def\Q{\mathbb Q_p}
\def\wh{\widehat}

\section{Introduction}
\label{Intro}

Let $C[0,1]$ be the set of continuous functions on $[0,1].$ Then the classical Bernstein polynomials of degree $n$
for $f\in C[0,1]$ are defined by
\begin{equation}\label{cla-def-ber}
\mathbb B_{n}(f)=\sum_{k=0}^nf\left(\frac kn\right)B_{k,n}(x),\quad 0\leq x\leq1
\end{equation}
where $\mathbb B_{n}(f)$ is called the Bernstein operator and
\begin{equation}\label{or-ber-poly}
B_{k,n}(x)=\binom nk x^k(x-1)^{n-k}
\end{equation}
are called the Bernstein basis polynomials (or the Bernstein polynomials of
degree $n$) (see \cite{SA}).
Recently, Acikgoz and Araci have studied the generating function for Bernstein
polynomials (see \cite{AA,AA2}). Their generating function for $B_{k,n}(x)$ is given by
\begin{equation}\label{AA-gen}
F^{(k)}(t,x)=\frac{t^ke^{(1-x)t}x^k}{k!}=\sum_{n=0}^\infty B_{k,n}(x)\frac{t^n}{n!},
\end{equation}
where $k=0,1,\ldots$ and $x\in[0,1].$ Note that
$$B_{k,n}(x)=\begin{cases}
\binom nk x^k(1-x)^{n-k}&\text{if } n\geq k \\
0,&\text{if }n<k
\end{cases}$$
for $n=0,1,\ldots$ (see \cite{AA,AA2}).

The Bernstein polynomials can also be defined in many
different ways. Thus, recently, many applications of these polynomials have been looked
for by many authors.
Some researchers have studied the Bernstein polynomials in the area of approximation theory (see \cite{AA,AA2,Be,KJY,Ph,SA}).
In recent years, Acikgoz and Araci \cite{AA,AA2} have introduced several type Bernstein polynomials.

In the present paper, we introduce the Bernstein polynomials on the ring of $p$-adic integers $\mathbb Z_p.$
We also investigate some interesting properties of the Bernstein polynomials related to the
bosonic $p$-adic integrals on the ring of $p$-adic integers $\mathbb Z_p.$

\section{Bernstein polynomials \\
related to the bosonic $p$-adic integrals on $\mathbb Z_p$}

Let $p$ be a fixed prime number.
Throughout this paper, $\mathbb Z_p, \mathbb Q_p$ and $\mathbb C_p$
will denote the ring of $p$-adic integers, the field of $p$-adic numbers and the completion
of the algebraic closure of $\mathbb Q_p,$ respectively.
Let $v_p$ be the normalized exponential valuation of $\mathbb C_p$ with
$|p|_p=p^{-1}.$
For $N\geq1,$ the bosonic distribution $\mu_1$ on $\mathbb Z_p$
\begin{equation}\label{mu}
\mu(a+p^N\Z)=\frac{1}{p^{N}}
\end{equation}
is known as the $p$-adic Haar distribution $\mu_{\text{Haar}},$ where $a+p^N\Z=\{ x\in\Q\mid |x-a|_p\leq p^{-N}\}$ (cf. \cite{KT1}).
We shall write $d\mu_{1}(x)$ to remind ourselves that $x$ is the variable
of integration.
Let $UD(\mathbb Z_p)$ be the space of uniformly differentiable function on $\mathbb Z_p.$
Then $\mu_{1}$ yields the fermionic $p$-adic $q$-integral of a function $f\in UD(\mathbb Z_p):$
\begin{equation}\label{Iqf}
I_{1}(f)=\int_{\Z} f(x)d\mu_{1}(x)=\lim_{N\rightarrow\infty}\frac{1}{p^{N}}
\sum_{x=0}^{p^N-1}f(x)
\end{equation}
(cf. \cite{KT1}).
Many interesting properties of (\ref{Iqf}) were studied by many authors
(cf. \cite{KT1,KCK} and the references given there).
For $n\in\mathbb N,$ write $f_n(x)=f(x+n).$ We have
\begin{equation}\label{de-3}
I_{1}(f_n)=I_{1}(f)+\sum_{l=0}^{n-1}f'(l).
\end{equation}
This identity is to derives interesting relationships involving Bernoulli numbers and polynomials.
Indeed, we note that
\begin{equation}\label{q-Euler-numb}
I_{1}((x+y)^n)=\int_{\Z}(x+y)^n d\mu_{1}(y)=B_n(x),
\end{equation}
where $B_n(x)$ are the Bernoulli polynomials (cf. \cite{KT1}).
From (\ref{or-ber-poly}), we have
\begin{equation}\label{b-ber-1}
\int_{\Z}B_{k,n}(x)d\mu_1(x)=\binom nk\sum_{j=0}^{n-k}\binom{n-k}j(-1)^{n-k-j}B_{n-j}
\end{equation}
and
\begin{equation}\label{b-ber-2}
\begin{aligned}
\int_{\Z}B_{k,n}(x)d\mu_1(x)&=\int_{\Z}B_{n-k,n}(1-x)d\mu_1(x) \\
&=\binom nk\sum_{j=0}^{k}\binom{k}j(-1)^{k-j}\sum_{l=0}^{n-j}\binom{n-j}{l}(-1)^lB_l.
\end{aligned}
\end{equation}
By (\ref{b-ber-1}) and (\ref{b-ber-2}), we obtain the following proposition.

\begin{proposition} For $n\geq k,$
$$\sum_{j=0}^{n-k}\binom{n-k}j(-1)^{n-k-j}B_{n-j}
=\sum_{j=0}^{k}\binom{k}j(-1)^{k-j}\sum_{l=0}^{n-j}\binom{n-j}{l}(-1)^lB_l.$$
\end{proposition}

From (\ref{q-Euler-numb}), we note that
\begin{equation}\label{b-ber-re}
B_n(2)=(B(1)+1)^n-n=(B+1)^n=B_n,\quad n>1
\end{equation}
with the usual convention of replacing $B^n$ by $B_n.$ Thus, we have
\begin{equation}\label{b-ber-re-2}
\begin{aligned}
\int_{\Z}x^nd\mu_1(x)&=\int_{\Z}(x+2)^nd\mu_1(x)-n \\
&=(-1)^n\int_{\Z}(x-1)^nd\mu_1(x)-n \\
&=\int_{\Z}(1-x)^nd\mu_1(x)-n
\end{aligned}
\end{equation}
for $n>1,$ since $(-1)^nB_n(x)=B_n(1-x).$
Therefore we obtain the following theorem.

\begin{theorem} For $n>1,$
$$\int_{\Z}(1-x)^nd\mu_1(x)=\int_{\Z}x^nd\mu_1(x)+n.$$
\end{theorem}

And also we obtain
\begin{equation}\label{b-ber-re-3}
\begin{aligned}
\int_{\Z}B_{n-k,k}(x)d\mu_1(x)&=\int_{\Z}x^{n-k}(1-x)^kd\mu_1(x) \\
&=\sum_{l=0}^{n-k}\binom{n-k}l(-1)^{l}\int_{\Z}(1-x)^{l+k}d\mu_1(x) \\
&=\sum_{l=0}^{n-k}\binom{n-k}l(-1)^{l}\left\{\int_{\Z}x^{l+k}d\mu_1(x)+l+k\right\} \\
&=\sum_{l=0}^{n-k}\binom{n-k}l(-1)^{l}(B_{l+k}+l+k).
\end{aligned}
\end{equation}
Therefore we obtain the following result.

\begin{corollary} For $k>1,$
$$\int_{\Z}B_{n-k,k}(x)d\mu_1(x)=\sum_{l=0}^{n-k}\binom{n-k}l(-1)^{l}(B_{l+k}+l+k).$$
\end{corollary}

From the property of the Bernstein polynomials of degree $n,$ we easily see that
\begin{equation}\label{b-ber-3}
\begin{aligned}
&\int_{\Z}B_{k,n}(x)B_{k,m}(x)d\mu_1(x) \\
&=\binom nk\binom mk\int_{\Z}x^{2k}(1-x)^{n+m-2k}d\mu_{1}(x)\\
&=\binom nk\binom mk
\sum_{l=0}^{n+m-2k}\binom{n+m-2k}l(-1)^{l}B_{2k+l}
\end{aligned}
\end{equation}
and
\begin{equation}\label{b-ber-4}
\begin{aligned}
&\int_{\Z}B_{k,n}(x)B_{k,m}(x)B_{k,s}(x)d\mu_1(x) \\
&=\binom nk\binom mk\binom sk\int_{\Z}x^{3k}(1-x)^{n+m-3k}d\mu_{1}(x)\\
&=\binom nk\binom mk\binom sk
\sum_{l=0}^{n+m+s-3k}\binom{n+m+s-3k}l(-1)^{l}B_{3k+l}.
\end{aligned}
\end{equation}
Continuing this process, we obtain the following theorem.

\begin{theorem}
The multiplication of the sequence of Bernstein polynomials
$$B_{k,n_1}(x),B_{k,n_2}(x),\ldots,B_{k,n_s}(x)$$
for $s\in\mathbb N$
with different degree under $p$-adic integral on $\mathbb Z_p$
can be given as
$$\begin{aligned}
&\int_{\Z}B_{k,n_1}(x)B_{k,n_2}(x)\cdots B_{k,n_s}(x)d\mu_1(x) \\
&=\binom{n_1}k\binom{n_2}k\cdots\binom{n_s}k
\sum_{l=0}^{n_1+n_2+\cdots+n_s-sk}\binom{n_1+n_2+\cdots+n_s-sk}l(-1)^{l}B_{sk+l}.
\end{aligned}$$
\end{theorem}

We put
$$B_{k,n}^{m}(x)=\underbrace{B_{k,n}(x)\times\cdots\times B_{k,n}(x)}_{m\text{-times}}.$$

\begin{theorem}
The multiplication of
$$B_{k,n_1}^{m_1}(x),B_{k,n_2}^{m_2}(x),\ldots,B_{k,n_s}^{m_s}(x)$$
Bernstein polynomials with different degrees $n_1,n_2,\cdots,n_s$ under $p$-adic integral on $\mathbb Z_p$
can be given as
$$\begin{aligned}
&\int_{\Z}B_{k,n_1}^{m_1}(x)B_{k,n_2}^{m_2}(x)\cdots B_{k,n_s}^{m_s}(x)d\mu_1(x) \\
&=\binom{n_1}k^{m_1}\binom{n_2}k^{m_2}\cdots\binom{n_s}k^{m_s}
\sum_{l=0}^{n_1m_1+n_2m_2+\cdots+n_sm_s-(m_1+\cdots+m_s)k}(-1)^{l} \\
&\quad\times\binom{n_1m_1+n_2m_2+\cdots+n_sm_s-(m_1+\cdots+m_s)k}lB_{(m_1+\cdots+m_s)k+l}.
\end{aligned}$$
\end{theorem}

\begin{theorem}
The multiplication of
$$B_{k_1,n_1}^{m_1}(x),B_{k_2,n_2}^{m_2}(x),\ldots,B_{k_s,n_s}^{m_s}(x)$$
Bernstein polynomials with different degrees $n_1,n_2,\cdots,n_s$
with different powers $m_1,m_2,\cdots,m_s$ under $p$-adic integral on $\mathbb Z_p$
can be given as
$$\begin{aligned}
&\int_{\Z}B_{k_1,n_1}^{m_1}(x)B_{k_2,n_2}^{m_2}(x)\cdots B_{k_s,n_s}^{m_s}(x)d\mu_1(x) \\
&=\binom{n_1}{k_1}^{m_1}\binom{n_2}{k_2}^{m_2}\cdots\binom{n_s}{k_s}^{m_s}
\sum_{l=0}^{n_1m_1+n_2m_2+\cdots+n_sm_s-(k_1m_1+\cdots+k_sm_s)}(-1)^{l} \\
&\quad\times\binom{n_1m_1+n_2m_2+\cdots+n_sm_s-(k_1m_1+\cdots+k_sm_s)}lB_{k_1m_1+\cdots+k_sm_s+l}.
\end{aligned}$$
\end{theorem}

\noindent
{\bf Problem.}
Find the Witt's formula for the Bernstein polynomials in $p$-adic number field.

\bibliography{central}

\end{document}